\documentclass[a4paper]{amsart}

\usepackage{amssymb,latexsym}
\usepackage{amsmath}
\usepackage{amsfonts}
\usepackage{amsthm}
\usepackage{amssymb}

\theoremstyle{plain}
\newtheorem{theorem}{Theorem}

\newtheorem{conjecture}{Conjecture}

\theoremstyle{definition}

\newtheorem{remark}{Remark}

\title[On the largest element in $D(n)$-quadruples]
{On the largest element in $D(n)$-quadruples}
\begin{document}

\date{}

%\footnotesize\date{\today}

\author[A. Dujella]{Andrej Dujella}
\address{
Department of Mathematics\\
Faculty of Science\\
University of Zagreb\\
Bijeni{\v c}ka cesta 30, 10000 Zagreb, Croatia
}
\email[A. Dujella]{duje@math.hr}

%\author[M. Kazalicki]{Matija Kazalicki}
%\address{
%Department of Mathematics\\
%Faculty of Science\\
%University of Zagreb\\
%Bijeni{\v c}ka cesta 30, 10000 Zagreb, Croatia
%}
%\email[M. Kazalicki]{matija.kazalicki@math.hr}

\author[V. Petri\v{c}evi\'c]{Vinko Petri\v{c}evi\'c}
\address{
Department of Mathematics\\
Faculty of Science\\
University of Zagreb\\
Bijeni{\v c}ka cesta 30, 10000 Zagreb, Croatia
}
\email[V. Petri\v{c}evi\'c]{vpetrice@math.hr}

\begin{abstract}
Let $n$ be a nonzero integer. A set of nonzero integers $\{a_1,\ldots,a_m\}$
such that  $a_ia_j+n$ is a perfect square for all $1\leq i<j\leq m$ is called
a $D(n)$-$m$-tuple. In this paper, we consider the question,
for given integer $n$ which is not a perfect square,
how large and how small can be the largest element in a $D(n)$-quadruple.
We construct families of $D(n)$-quadruples in which the largest element is of order
of magnitude $|n|^3$, resp. $|n|^{2/5}$.
\end{abstract}


\maketitle

\section{Introduction}
For a nonzero integer $n$, a set of distinct nonzero integers $\{a_1,a_2,\ldots,a_m\}$
such that  $a_ia_j+n$ is a perfect square for all $1\leq i<j\leq m$ is called
a $D(n)$-$m$-tuple (or a Diophantine $m$-tuple with the property $D(n)$).

The most studied case is $n=1$ and $D(1)$-$m$-tuples are called Diophantine $m$-tuples.
Fermat found the first Diophantine quadruple, it was the set $\{1,3,8,120\}$.
In 1969, Baker and Davenport \cite{BD69} proved that the set $\{1,3,8\}$ can be extended to a Diophantine quintuple only by adding $120$ to the set.
In 2004, Dujella~\cite{d05} proved that there are no Diophantine sextuples and that there are at most finitely many Diophantine quintuples.  Recently, He, Togb\'e and Ziegler proved that there are no Diophantine quintuples~\cite{HTZ16+} (see also \cite{BTF}). 
On the other hand, there are examples
of $D(n)$-quintuples and sextuples for $n\neq 1$, 
e.g. $\{8, 32, 77, 203, 528\}$ is a $D(-255)$-quintuple \cite{duje-acta2}, 
while $\{99, 315, 9920, 32768, 44460, 19534284\}$ is a $D(2985984)$-sextuple \cite{Gibbs1} 
(see also \cite{DKMS}). 
For an overview of results on Diophantine $m$-tuples and its generalizations see \cite{Duje-Notices}.

Several authors considered the problem of the existence of Diophantine quadruples
with the property $D(n)$. It is easy to show that there are no $D(n)$-quadruples
if $n\equiv 2 \pmod{4}$ (\cite{Brown,GS,MR}). Indeed, assume that
$\{a_1, a_2, a_3, a_4\}$ is a $D(n)$-quadruple. Since the square of an integer is
$\equiv  0$ or $1 \pmod{4}$, we have that
$a_i a_j \equiv 2$ or $3 \pmod{4}$. This implies that none of the
$a_i$'s is divisible by $4$. Therefore, we may assume that $a_1 \equiv a_2 \pmod{4}$.
But now we have that $a_1 a_2 \equiv 0$ or $1 \pmod{4}$, a contradiction.
On the other hand, it is shown in \cite{duje-acta1} that
if $n\not\equiv 2\pmod{4}$ and
$n\not\in S =\{-4, -3, -1, 3, 5, 8, 12, 20\}$, then there exists at least one $D(n)$-quadruple.
For $n\in S$, the question of the existence of $D(n)$-quadruples is still open. 

The Lang conjecture on varieties of general type implies 
that the size of sets with the property $D(n)$ is bounded by an absolute constant
(independent on $n$). It is known that the size of sets with the property $D(n)$
is $\leq 31$ for $|n|\leq 400$; $< 15.476 \log|n|$ for $|n| > 400$, and $ < 3 \cdot 2^{168}$ for $n$ prime (see \cite{Duje-size1,Duje-size2,DL-IMRN} and also \cite{murty}). 

It is easy to see that there exist infinitely many $D(1)$-quadruples.
Indeed, the set $\{k-1,k+1,4k,16k^3-4k\}$ for $k\geq 2$ is a $D(1)$-quadruple
(see e.g. \cite{duje-pdeb1}). More precisely, it was proved in \cite{M-S} that
the number of $D(1)$-quadruples with elements $\leq N$ is $\sim C\sqrt[3]{N}\log N$,
where $C\approx 0.338285$ (the main contribution comes from the quadruples of the
form $\{a,b,a+b+2r,4r(a+r)(b+r)\}$ where $ab+1=r^2$, see \cite{duje-raman}).
If $n$ is a perfect square, say $n=\ell^2$, then by multiplying
elements of a $D(1)$-quadruple by $\ell$ we obtain a $D(\ell^2)$-quadruple,
and thus we conclude that there exist infinitely many $D(\ell^2)$-quadruples.
Moreover, it was proved in \cite{duje-acta1} that any $D(\ell^2)$-pair $\{a,b\}$,
such that $ab$ is not a perfect square, can be extended to a $D(\ell^2)$-quadruple.

The following conjecture was proposed in \cite{duje-darf}.

\begin{conjecture} \label{conj1}
If a nonzero integer $n$ is not a perfect square,
then there exist only finitely many $D(n)$-quadruples.
\end{conjecture}

As we already mentioned, it is easy to verify the conjecture
in case $n\equiv 2\pmod{4}$ since there does not exist a $D(n)$-quadruple in that case.
Only other cases where the conjecture is known to be true are the cases $n=-1$ and $n=-4$, see \cite{DFF,DF-London}
(these two cases are equivalent since, by \cite{duje-acta1}, all elements of a $D(-4)$-quadruple
are even).

Motivated by Conjecture \ref{conj1}, in this paper we consider the question,
for given integer $n$ which is not a perfect square,
what can be said about the largest element in a $D(n)$-quadruple.
In particular, the question how large it can be (compared with $|n|$)
is closely related with Conjecture \ref{conj1}.
On the other hand, the question how small it can be (again compared with $|n|$)
makes sense also in the case when $n$ is a perfect square.
Since $\{a,b,c,d\}$ is a $D(n)$-quadruple if and only if $\{-a,-b,-c,-d\}$ has the same property,
without loss of generality we may assume that $\max\{|a|,|b|,|c|,|d|\}=d$.
Our main results are collected in the following theorem.

\begin{theorem} \label{tm1}
Let $\delta,\varepsilon$ be real numbers such that
$2/5 \leq \delta \leq 3$ and $\varepsilon >0$. Then there exist an integer $n$
which is not a perfect square and a $D(n)$-quadruple
$\{a,b,c,d\}$ such that
$$ \left|\frac{\log(\max\{|a|,|b|,|c|,|d|\})}{\log{|n|}} - \delta \right| < \varepsilon. $$
\end{theorem}

In Section \ref{sec2} we will consider $D(n)$-quadruples, where $n$ is not a perfect square,
with large elements
and construct family of quadruples with $d$ of order of magnitude $|n|^3$,
while in Section \ref{sec3} we will consider $D(n)$-quadruples with small elements
and construct family of quadruples with $d$ of order of magnitude $|n|^{2/5}$.
Since elements of both families of quadruples are polynomials in one variable,
a standard construction with $D(n)$-quadruples will finish the proof of Theorem \ref{tm1}.

It should be noted that we do not know what are best possible results in both direction,
i.e. is there any family of $D(n)$-quadruples with $d$ of order of magnitude $|n|^{\delta}$
with $\delta > 3$ or $\delta < 2/5$. We will show in Section \ref{sec3} that we cannot have $\delta < 1/4$.

\section{$D(n)$-quadruples with large elements} \label{sec2}

The proof of the fact that for $n\not\equiv 2\pmod{4}$ and
$n\not\in \{-4, -3, -1, 3, 5, 8, 12, 20\}$ there exists at least one $D(n)$-quadruple \cite{duje-acta1},
is based on explicit formulas for $D(n)$-quadruples, where
$n=4k+3$, $n=8k+1$, $n=8k+5$, $n=8k$, $n=16k+4$ and $n=16k+12$, while elements of
$D(n)$-quadruples are polynomials in $k$. For example,
$$ \{1, 9k^2 + 8k + 1, 9k^2 + 14k + 6, 36k^2 + 44k + 13\} $$
is a $D(4k+3)$-quadruple.
This example shows that it is possible to have $\max\{a,b,c,d\}$ $\sim \frac{9}{4} |n|^2$.
The same conclusion also follows from the fact that
\begin{equation} \label{eq:duje-ana}
\begin{gathered}
\{1, 144k^4 + 216k^3 + 113k^2 + 20k + 1, 144k^4 + 360k^3 + 329k^2 + 134k + 22, \\
 576k^4 + 1152k^3 + 848k^2 + 272k + 33\}
\end{gathered}
\end{equation}
is a $D((4k + 1)(4k + 3))$-quadruple
 (see \cite{duje-ana}).

\medskip

In order to find families of quadruples with larger elements (compared with $|n|$),
we performed an extensive search for $D(n)$-quadruples (similar as explained in \cite{duje-vinko-dn1dn2})
and then sieve them according to the requirement that
$\max\{a,b,c,d\}/|n|^2$ is relatively large (in particular, larger that $9/4$).
Then we search for properties which are common to several found quadruples.
For example, quadruples containing one small element (e.g. $a=1$) and quadruples
containing a regular triple (a $D(n)$-triple $\{b,c,d\}$ is called regular if
$(b+c-d)^2=4(bc+n)$). We have extracted the following interesting examples:

\begin{center}
\begin{tabular}{c|c|c}
	$\{a,b,c,d\}$ & $n$ 		& $d/n^2$ \\
	\hline
$\{1, 2912, 131977, 174097\}$  &   $-208$   & $4.024062$ \\
$\{1, 16896, 1980161, 2362881\}$ & $-512$   & $9.013676$ \\
$\{1, 56640, 12525465, 14266673\}$ & $-944$   &$16.009535$ \\
	\end{tabular}
	\end{center}

\medskip

These examples suggest that for every positive integer $k$ there might exist
a quadruple such that $d/n^2\approx k^2$. Furthermore, in these three examples
(and some other found examples) we have that $n|b$. In fact, these examples suggests
that we may take that
$b/n=-k(4k-1)$. Starting with these assumptions, it is now easy to reconstruct the
corresponding family $\{a,b,c,d\}$ of $D(n)$-quadruples:
\begin{align*}
a &=1, \\
b &=256k^4-128k^3-48k^2+16k, \\
c &=4096k^6-4096k^5-512k^4+1152k^3+16k^2-88k-7, \\
d &=4096k^6-2048k^5-1792k^4+640k^3+288k^2-48k-15, \\
n &=-64k^2+16k+16.
\end{align*}
Here $k$ is an arbitrary nonzero integer. By taking $|k|\to \infty$, we obtain
$d/n^2 \to \infty$ and $d/|n|^3 \to 1/64$.

After the substitution $k=z+1/8$, the family of quadruples becomes
\begin{equation*}
\begin{gathered}
 \{1, 256z^4-72z^2+\frac{17}{16},
 4096z^6-1024z^5-2112z^4+416z^3+335z^2-\frac{153}{4}z-\frac{1007}{64}, \\
 4096z^6+1024z^5-2112z^4-416z^3+335z^2+\frac{153}{4}z-\frac{1007}{64}\},
\end{gathered}
\end{equation*}
 with $n=-64z^2+17$. Thus, $n$, $b$ and $d+c$ are even polynomials in $z$, while $d-c$ is odd
 (analogously as in (\ref{eq:duje-ana})).
From these properties, it is also possible to reconstruct the polynomials by the method of
undetermined coefficients.

\medskip

In our example we have $n \sim -64k^2$ and $d \sim 4096 k^6$. Hence, $\log{d}/\log{|n|} \to 3$
as $k\to \infty$.
If we take $k=y^{\ell_1}$, and multiply all elements of the quadruple by $y^{\ell_2}$,
we get quadruple in which $n \sim -64 y^{2\ell_1+2\ell_2}$ and $d \sim 4096 y^{6\ell_1+\ell_2}$.
Hence, now we have $\log{d}/\log{|n|} \sim (6\ell_1+\ell_2)/(2\ell_1+2\ell_2)$.
By varying nonnegative integers $\ell_1$ and $\ell_2$, we get that any point from the interval
$[1/2, 3]$ is an accumulation point of the set
$\{\log{d}/\log{|n|} \,:\, \mbox{$\{a,b,c,d\}$ is a $D(n)$-quadruple for some non-square $n$}\}$.

\section{$D(n)$-quadruples with small elements} \label{sec3}

In this section, we consider the question how small can be elements of a $D(n)$-quadruple,
in particular how small can be its largest element.
As we have already seen at the end of the previous section,
it is easy to get quadruples in which $\max\{a,b,c,d\}$ is of order of magnitude $|n|^{1/2}$.
Indeed, we can take any nonzero integer $k$ and any
fixed $D(k)$-quadruple $\{a_1,b_1,c_1,d_1\}$,
%$|a_1|<|b_1|<|c_1|<|d_1|$,
 and multiply its elements by
large positive integer $\ell$ to get $D(k \ell^2)$-quadruple $\{a_1 \ell, b_1 \ell, c_1 \ell, d_1 \ell\}$,
%which yields $d/\sqrt{|n|} = d_1 \ell / \sqrt{|k \ell^2|} \to d_1/\sqrt{|k|}$ as $\to \infty$.
which yields $\log(\max\{|a|,|b|,|c|,|d|\})/\log{|n|} \to 1/2$ as $\ell \to \infty$.

Thus, we are interested if there are (families of) examples with $d$ significantly smaller that $|n|^{1/2}$.
If $n<0$, then we can assume that $0<a<b<c<d$,
and from $cd+n>0$ it follows that $d>|n|^{1/2}$.
Hence, we may assume that $n>0$. We claim that $d$ cannot be smaller that $n^{1/4}$.
Indeed, let $|a|\leq |b|\leq |c|\leq d < n^{1/4}$. Since $d\geq 2$, we have $n>16$.
We may assume that $c>b$ (if $c<b$ the proof is analogous).
From $cd+n=r^2$, $bd+n=s^2$ we get $n^{1/2} > d(c-b) = r^2-s^2 \geq (s+1)^2 -s^2 = 2s+1$.
We obtain $s< \frac{1}{2}n^{1/2}$, which implies $bd < -\frac{3}{4}n$, and this contradicts $|bd| < n^{1/2}$.

As the examples of $D(16k^4-72k^2+48k+9)$-triple $\{4k-4, 8k-4, 12k\}$ and
$D(144k^4+264k^3+181k^2+52k+5)$-triple $\{-6k-1, 2k+1, 6k+4\}$ show,
in a $D(n)$-triple we may have all elements of the order $n^{1/4}$.
However, we were not able to find $D(n)$-quadruples with the same property.

\subsection{$D(n)$-quadruples of the form $\{a,-a,b,-b\}$}

In considering certain problems with $D(n)$-quadruples, it might be convenient to
study sets of the form $\{a,-a,b,-b\}$. In order to satisfy the definition of a $D(n)$-quadruple,
such sets has to satisfy only four conditions: $-a^2+n$, $-b^2+n$, $ab+n$ and $-ab+n$ are
perfect squares, compared with six conditions which has to be satisfied by general set of four elements.

By considering $D(n)$-quadruples of the form $\{a,-a,b,-b\}$, we found examples
with $\max\{|a|,|b|\} \leq n^{1/2}$:
$$ \{-4u, 4u, -3-u^2,  3+u^2\} $$
is a $D((u^2+9)(1+u^2))$-quadruple,
$$ \{-4u(u-1)(u-2), 4u(u-1)(u-2), -(u^2-2u+2)^2, (u^2-2u+2)^2\} $$
is a $D((u^2-2u+2)^4)$-quadruple, while
$$ \{-4u^2-2u-1, 4u^2+2u+1, -4u(u+1), 4u(u+1)\} $$
is a $D((10u^2+2u+1)(2u^2+2u+1))$-quadruple.

We can improve slightly these results to get families of
$D(n)$-quadruples of the form $\{a,-a,b,-b\}$ such that $\max\{|a|,|b|\}$
is of order of magnitude $n^{9/20}$.

Let $-a^2+n=r^2$ and $ -ab+n=(r-1)^2$. We get
$$ n = a^2+r^2, \quad r = (ab-a^2+1)/2. $$
We write the third condition in the form
$-b^2+n=(r-t)^2$, which gives that
$a^2 t^2+4t+4a^2-4a^2t-4t^2$ is a perfect square, say
$$ (t^2+4-4t)a^2-4t^2+4t = ((t-2)a+u)^2. $$
We get
$$ a = -(4t^2-4t+u^2)/(2u(t-2)), \quad b = -(4t-8t^2+4t^3+u^2)/(2u(t-2)). $$
It remains the satisfy the last condition that $ab+n$ is a perfect square. The condition leads to
\begin{gather*}
 16 t^8-64 t^7+8 u^2 t^6+96 t^6-64 t^5+16 t^4-40 u^2 t^4 +u^4 t^4+48 u^2 t^3+4 u^4 t^3 \\
\mbox{}-16 u^2 t^2 -8 u^4 t+4 u^4+2 u^6=\Box.
\end{gather*}
If we take $u=t(t-1)/v$,
the condition becomes
$$ (v^2+2)t^4+(4v^2-4)t^3+(8v^4+2)t^2+(16v^4-8v^2)t+16v^6-16v^4+4v^2=s^2. $$
This quartic over $\mathbb{Q}(v)$ has a $\mathbb{Q}(v)$-rational point
$P_1=[0, -2v(2v^2-1)]$, and therefore it can be in standard way (see e.g. \cite{connell})
transformed into an elliptic curve. There is another point on the quartic:
$$ P_2=[-4(v-1)v/(2v+1), -2v(16v^4-16v^3+14v^2-8v+3)/((2v+1)^2)]. $$
If we take the point $2P_2$, we get
$$ t=\frac{-8v(16v^4-16v^3-2v^2+8v-3)}{32v^3-56v^2+20v+1}, $$
which gives the $D(n)$-quadruple $\{a,-a,b,-b\}$, where
\begin{align*}
a &= 16384 v^9-32768 v^8+20480 v^7-5120 v^6+4608 v^5-4096 v^4+1216 v^3-304 v^2 \\
&\,\,\,\,\mbox{}+164 v-24, \\
b &= 12288 v^8-24576 v^7+21504 v^6-11520 v^5+2880 v^4+576 v^3-528 v^2+156 v-24, \\
n & =268435456 v^{20}-1073741824 v^{19}+1879048192 v^{18}-1946157056 v^{17} \\
&\,\,\,\,\mbox{}+1392508928 v^{16}-788529152 v^{15}+465567744 v^{14}-412090368 v^{13} \\
&\,\,\,\,\mbox{}+412483584 v^{12}-328990720 v^{11}+205324288 v^{10}-110215168 v^9 \\
&\,\,\,\,\mbox{}+53587968 v^8-22474752 v^7+7394304 v^6-1852160 v^5+468752 v^4 \\
&\,\,\,\,\mbox{}-164480 v^3+47872 v^2-7552 v+580,
\end{align*}
which satisfies $\log (\max\{|a|,|b|\})/\log{n} \to 9/20$ as $v \to \infty$.

\subsection{$D(n)$-quadruples with $d^5 \sim 8n^2$}

Now we describe the construction of a family of $D(n)$-quadruples in which $d$ is of
order of magnitude $n^{2/5}$. The construction is motivated by the following experimentally found example:
$$ \{468,335,-85,-448\} $$
is a $D(1312164)$-quadruple.

Let $n=x^2+y$ and put $ab+n=(x+s)^2$, $bd+n=(x-s)^2$. We get
$$y = -ab+2xs+s^2, \quad d = \frac{ab-4xs}{b}. $$
Now put $ac+n=(x-r)^2$, $cd+n=(x+r)^2$, and we get
$$ a = \frac{s(s+2x-r)}{b}, \quad c = \frac{-br}{s}. $$
By taking $bc+n=(x-t)^2$, we obtain
$$ x = \frac{b^2 r-s^2 r+s t^2}{2st}. $$
It remains to satisfy the condition that $ad+n$ is a perfect square.
We put $ad+n=(x-r-s-t+1)^2$. This leads to the equation
\begin{equation} \label{eq:bs}
\begin{gathered}
-2b^2s^2t^2-b^2t^3s+b^2st^2-2s^3r^2b^2+sb^4r^2-s^5t^2+s^3b^2t^2+b^2s^2t^3-ts^2b^2r \\
\mbox{}+s^3t^4+tb^4r+s^5r^2+4b^2rs^2t^2-b^4rt^2-tb^4rs+ts^3b^2r-s^3r^2t^2-b^4r^2t \\
\mbox{}+b^2s^2r^2t+b^2st^3r+b^2r^2st^2-2b^2rst^2=0,
\end{gathered}
\end{equation}
In order to find some of its solution,
we introduce the condition $x-y=t/2$.
This gives $t =  \frac{-s^2+b^2}{2s^2}$, and by inserting it in (\ref{eq:bs})
(the discriminant of the equation in $r$ becomes a square), we get
% (b^4+2*b^2*s^3+2*b^2*r*s^2-4*s^2*b^2-4*s^5*r-2*s^5-2*s^4*r-s^4)=0
$$ r = \frac{-(b^4+2b^2s^3-4s^2b^2-2s^5-s^4)}{2s^2(-2s^3-s^2+b^2)}. $$
Thus, we obtain the rational $D(n)$-quadruple $\{a,b,c,d\}$, where
\begin{align*}
a &= \frac{-b(s-1)(2s^3-3s^2+b^2)}{s(-2s^3-s^2+b^2)}, \\
c &= \frac{b(b^4+2b^2s^3-4s^2b^2-2s^5-s^4)}{2s^3(-2s^3-s^2+b^2)}, \\
d &= \frac{(b^2+2s^3+s^2)(b^2-s-2s^2)}{b(-2s^3-s^2+b^2)}, \\
n &= (b+s)(b-s)(2b^3s-b^3+3b^2s-2s^2b^2-3s^2b-4bs^3+s^3+4bs^4-4s^5) \\
\,\,\,\,\,\,\,\,&\mbox{}\times(2b^3s-b^3-3b^2s+2s^2b^2-3s^2b-4bs^3-s^3+4bs^4+4s^5) \\
\,\,\,\,\,\,\,\,&\mbox{}\times(16s^4(-2s^3-s^2+b^2)^2)^{-1}.
\end{align*}

%It remains to choose the parameters $b$ and $s$ in such a way that the ratio
%$\log (\max\{|a|,|b|,|c|,|d|\}) / \log{n}$ becomes relatively small.

In order to get quadruples with integers elements, we put $b=ks$ and search for the solution in the form
$s=s_3 k^3+s_2 k^2+s_1k+s_0$. The condition that rational functions appearing in $a,c,d,n$
become polynomials, leads to $s_0 = -1/2$, $s_2 = 1/2$, $s_3 = -s_1/3$.
We take $s_1=-3/2$ and $k=2v-1$, and we obtain the
$D(64v^{10}-128v^9-64v^8+240v^7-32v^6-136v^5+41v^4+22v^3-7v^2])$-quadruple
$$  \{8v^4-4v^3-8v^2, 8v^4-12v^3+4v-1, -4v^3+2v^2+2v-1, -8v^4+4v^3+12v^2-4v-4 \}, $$
which satisfies $\log (\max\{|a|,|b|,|c|,|d|\})/\log{n} \to 2/5$ as $v \to \infty$.
For $v=3$ we get our starting motivating example: $D(1312164)$-quadruple $\{468,335,-85,-448\}$.

\medskip

We now apply the same argument as at the end of Section \ref{sec2}.
We have $n \sim 64v^{10}$ and $\max\{|a|,|b|,|c|,|d|\} \sim 8v^4$.
If we take $v=y^{\ell_1}$, and multiply all elements of the quadruple by $y^{\ell_2}$,
we get quadruple in which $n \sim 64 y^{10\ell_1+2\ell_2}$ and
$\max\{|a|,|b|,|c|,|d|\} \sim 8 y^{4\ell_1+\ell_2}$.
Hence, we have $\log(\max\{|a|,|b|,|c|,|d|\})/\log{|n|} \sim (4\ell_1+\ell_2)/(10\ell_1+2\ell_2)$.
By varying nonnegative integers $\ell_1$ and $\ell_2$, we get that any point from the interval
$[2/5, 1/2]$ is an accumulation point of the set
$$ \Big\{\frac{\log(\max\{|a|,|b|,|c|,|d|\})}{\log{|n|}} \,:\, \mbox{$\{a,b,c,d\}$ is a $D(n)$-quadruple for some non-square $n$} \Big\}, $$
which together with the mentioned result from the end of Section \ref{sec2}
finishes the proof of Theorem \ref{tm1}.

\begin{remark} \label{rem1}
As we mentioned in the introduction, the question how small the largest element of a
$D(n)$-quadruple can be make sense also in the case when $n$ is a perfect square.
By taking the $D(42849)$-quadruple $\{188, 140, -160, -198\}$ as a motivating example,
we found that
$$ \{60 u^2-24 u-4, 100 u-60, -4 (5 u-2) (3 u-1), -90 u^2+96 u-30\} $$
is a $D((15 u-7)^2 (5 u-1)^2)$-quadruple,
and this yields $\log (\max\{|a|,|b|,|c|,|d|\})/\log{n} \to 1/2$ as $u \to \infty$.
Thus, a version of Theorem \ref{tm1} with $n$ a nonzero perfect square holds for
$1/2 \leq \delta < \infty$. Again, we do not know whether the lower bound for $\delta$ is best possible.
\end{remark}

\bigskip

{\bf Acknowledgements.}
The authors were supported by the Croatian Science Foundation under the project no.~IP-2018-01-1313.
The authors acknowledge support from the QuantiXLie Center of Excellence, a project
co-financed by the Croatian Government and European Union through the
European Regional Development Fund - the Competitiveness and Cohesion
Operational Programme (Grant KK.01.1.1.01.0004).
The authors acknowledge the usage of the supercomputing resources
of Division of Theoretical Physics at Ru\dj{}er Bo\v{s}kovi\'c Institute.

\end{document}